\RequirePackage{amsthm}
\documentclass[sn-mathphys,Numbered]{sn-jnl}

\usepackage{graphicx}%
\usepackage{multirow}%
\usepackage{amsmath,amssymb,amsfonts,mathtools}%
\usepackage{mathrsfs}%
\usepackage[title]{appendix}%
\usepackage{xcolor}%
\usepackage{textcomp}%
\usepackage{manyfoot}%
\usepackage{booktabs}%
\usepackage{algorithm}%
\usepackage{algorithmicx}%
\usepackage{algpseudocode}%
\usepackage{listings}%

\theoremstyle{thmstyleone}%
\newtheorem{theorem}{Theorem}
\theoremstyle{thmstyletwo}%
\newtheorem{fact}{Fact}

\theoremstyle{thmstylethree}%
\newtheorem{definition}{Definition}%

\begin{document}

\title[]{Two-color exchangeability does not imply exchangeability}


\author*[]{\fnm{Raghavendra} \sur{Tripathi}}\email{raghavt@uw.edu}



\affil*[]{\orgdiv{Department of Mathematics}, \orgname{University of Washington},  \city{Seattle}, \postcode{}, \state{Washington}, \country{USA}}




\abstract{We provide simple examples of two-color exchangeable sequences $\xi=(\xi_1, \xi_2, \ldots, \xi_n)$ that are not exchangeable. This answers a question of Bladt and Shaiderman~\cite[Question 2.6]{bladt2019characterisation} for finite two-color exchangeable sequences.}




\maketitle
\section{Introduction}
A finite or infinite sequence of random variables $(\xi_i)_{i\geq 1}$ taking values in a Polish space is \emph{exchangeable} if \[(\xi_1, \ldots, \xi_n)\stackrel{d}{=}(\xi_{\pi(1)}, \ldots, \xi_{\pi(n)})\]
for every $n$ and every permutation $\pi$ on $[n]$. Exchangeability is possibly one of the most widely studied topic in probability since the work of de Finetti's in 1930. Given a finite or infinite sequence of random variables $\xi_i, i=1, 2, \ldots$, one can associate a sequence of random empirical measures as $\mu_n\coloneqq n^{-1}\sum_{i=1}^{n}\delta_{\xi}$. The following theorem stated in~\cite{kallenberg2005probabilistic} characterizes exchangeable sequences in terms of associated sequence of random empirical. 
\begin{theorem}[Kallenberg]
A sequence $\xi\coloneqq (\xi_n)_{n\geq 1}$ is exchangeable if and only if the associated sequence of empirical measures $\mu\coloneqq (\mu_n)_{n\geq 1}$ forms a reverse martingale (with respect to the filtration $\mathcal{T}_n=\sigma(\mu_n, \mu_{n+1}, \ldots)$).
\end{theorem}
Although the above result is correct, the proof presented in~\cite{kallenberg2005probabilistic} was not completely rigorous. The first complete proof of the above result was presented only very recently in~\cite{bladt2019characterisation}. The authors in~\cite{bladt2019characterisation} introduce a notion of \emph{two-color exchangeability} which seems interesting in its own right. We present the definition below~\cite[Definition 2.4]{bladt2019characterisation} for completeness. We also refer the reader to~\cite[Section 2.2]{bladt2019characterisation} for an interesting connection of two-color exchangeability with Strassen's marginal problem~\cite{strassen1965existence, edwards1978existence}.

\begin{definition}[Two-color exchangeability]
    A finite or infinite sequence of random variables $\xi\coloneqq (\xi_i)_{i\geq 1}$ taking values in a Polish space $(S, \mathcal{S})$ is said to be two-coloring exchangeable if $(f(\xi_i))_{i\geq 1}$ is a binary exchangeable sequence for every measurable function $f:S\to \{0, 1\}$.
\end{definition}

It is easy to see that if a sequence $\xi\coloneqq (\xi_i)_{i\geq 1}$ is exchangeable then it is also two-color exchangeable. The author in~\cite[Question 2.6]{bladt2019characterisation} ask if the converse is true or not? 

The purpose of this short note is to answer the above question in negative for finite sequences. In particular, for each $n\geq 2$, we construct a sequence of random variables $\xi_1, \xi_2, \ldots, \xi_n$ taking values in $\{-1, 0, 1\}$ such that $\xi\coloneqq (\xi_1, \xi_2, \ldots, \xi_n)$ is two-color exchangeable but not exchangeable. For pedagogical reasons, we first give such an example for $n=2$ case in Section~\ref{sec:n=2} and then in Section~\ref{sec:General} we construct the example for general $n\geq 3$.

\section{Case \texorpdfstring{$n=2$}{n=2}}\label{sec:n=2}
Let $\xi=(\xi_1, \xi_2)$ where $\xi_1, \xi_2$ are random variables taking values in $\{-1, 0, 1\}$. We now specify the distribution of $\xi$. Suppose,
\begin{align*}
    \mathbb{P}(\xi=(x, x))&=0,\quad \forall x\in \{-1, 0, 1\}\\
    \mathbb{P}(\xi=(-1, 0))&=\mathbb{P}(\xi=(0, 1))=\mathbb{P}(\xi=(1, -1))=\frac{1}{9}\\
    \mathbb{P}(\xi=(0, -1))&=\mathbb{P}(\xi=(1, 0))=\mathbb{P}(\xi=(-1, 1))=\frac{2}{9}\;.\\
\end{align*}

Note that $\xi$ is {\bf not} exchangeable since 
\[\mathbb{P}(\xi=(-1, 0))=\mathbb{P}((\xi_1, \xi_2)=(-1, 0))\neq \mathbb{P}((\xi_2, \xi_1)=(-1, 0))=\mathbb{P}(\xi=(0, -1))\;.\]

However, we will show that $\xi$ is two-color exchangeable. To this end, we need to check that $(f(\xi_1), f(\xi_2))$ is exchangeable for every $f:\{-1, 0, 1\}\to \{0, 1\}$. We begin by noting that there are a total of $8$ such function. For $f\equiv 0$ or $f\equiv 1$, we trivially have that $(f(\xi_1), f(\xi_2))$ is exchangeable. 

For $i\in \{-1, 0, 1\}$, define $f_i:\{-1, 0, 1\}\to \{0, 1\}$ by $f_i=\delta_i$ where $\delta_i(x)=1$ if $x=i$ and $0$ otherwise. And, define $g_i = 1-f_i$. It is easy to check that these are all (non-trivial) binary functions on $\{-1, 0, 1\}$. It is also easy to see that $(f_i(\xi_1), f_i(\xi_2))$ is exchangeable if and only if $(g_i(\xi_1), g_i(\xi_2))$ is exchangeable. Therefore, it suffices to check that $(f_i(\xi_1), f_i(\xi_2))$ is exchangeable for $i\in \{-1, 0, 1\}$. This is done by simple calculation that we furnish below for completeness. 

In the following, we write $\mu_{i}$ for the probability mass function of $(f_i(\xi_1), f_i(\xi_2))$ for simplicity. Note that 
\begin{align*}
    \mu_{-1}((1, 1)) &= \mathbb{P}(\xi=(-1, -1)) = 0\\
    \mu_{-1}((1, 0)) &= \mathbb{P}(\xi=(-1, 0))+ \mathbb{P}(\xi=(-1, 1)) = \frac{1}{9}+\frac{2}{9} = \frac{1}{3}\\
     \mu_{-1}((0, 1)) &= \mathbb{P}(\xi=(0, -1))+ \mathbb{P}(\xi=(1, -1)) = \frac{2}{9}+\frac{1}{9} = \frac{1}{3}\\
     \mu_{-1}((0, 0)) &= \frac{1}{3}\;.
\end{align*}

Since $\mu_{-1}((1, 0))=\mu_{-1}((0, 1))$, it follows that $(f_{-1}(\xi_1), f_{-1}(\xi_2))$ is exchangeable. Similar calculations hold for $f_{0}$ and $f_{1}$ as well. This completes the proof.

\subsection{A word about the above counterexample}
It might be instructive to say a few words about how we arrived at the above counterexample. To this end, we started with the simplest case of a sequence of two random variables $\xi_1, \xi_2$ taking values in a discrete set with three elements. For simplicity, we take the range of the random variables to be $\{-1, 0, 1\}$. 

Let $\xi=(\xi_1, \xi_2)$. Since $\xi$ can take only $9$ distinct values, the probability mass function of $\xi$ is completely determined by an element $\mathbf{p}\in \Delta_{9}$ where $\Delta_{9}$ is the unit simplex $\Delta_{9}=\{(p_1, \ldots, p_9): p_1+\ldots p_9 =1, p_i\geq 0\}$. Let $\mu$ denote the probability mass function of $\xi$ and denote 
\begin{align*}
    x_1 &= \mu((-1, -1)) &&\\
    x_2 &= \mu((-1, 0)), &&\widetilde{x}_2 = \mu((0, -1))\\
    x_3 &= \mu((0, 0)), \\
    x_4 &= \mu((0, 1)), &&\widetilde{x}_4 = \mu((1, 0))\\
    x_5 &= \mu((1, 1))\\
     x_6 &= \mu((1, -1)), &&\widetilde{x}_6 = \mu((-1, 1))\;.
\end{align*}

Any arbitrary choice of non-negative real numbers $x_1, \ldots, x_6$ and $\widetilde{x}_2, \widetilde{x}_4, \widetilde{x}_6$ such that $x_1+\ldots+x_6+\widetilde{x}_2+\widetilde{x}_4+\widetilde{x}_6=1$ yields a unique probability mass function and hence uniquely determines a random vector $\xi=(\xi_1, \xi_2)$. The following fact is straightforward.
\begin{fact}
    Let $\xi$ be as above. Then, $\xi$ is exchangeable if and only if $x_j=\widetilde{x}_j$ for $j=2, 4, 6$. 
\end{fact}

On the other hand, doing a calculation similar to the one done in the above subsection yields the following. 

\begin{fact}
   $\xi$ is two-color exchangeable if and only if 
   \begin{align*}
       x_2+\widetilde{x}_6 &= \widetilde{x}_2+ x_6\\
       x_2+\widetilde{x}_4 &= \widetilde{x}_2+ x_4\\
       x_4+\widetilde{x}_6 &= \widetilde{x}_4+ x_6\;.
   \end{align*}
\end{fact}

Now it is easy to construct examples of two-color exchangeable sequences of length $2$ that are not exchangeable.

\section{Case \texorpdfstring{$n\geq 3$}{n>=3}}\label{sec:General}
Throughout this section, we fix $n\geq 3$. We will construct $\xi=(\xi_1, \ldots, \xi_n)$ such that each $\xi_i$ takes values in $\{-1, 0, 1\}$ such that $\xi$ is two-color exchangeable but not exchangeable. Note that $\xi$ takes values in $\{-1, 0, 1\}^{n}$. However, for simplicity we will only consider the random vectors that take values in $\Omega$ where 
 \[\Omega = \{(x_1, \ldots, x_n)\in \{-1, 0, 1\}^{n}: \text{Each symbol appears exactly } 0\text{ or} 1\text{ or }n-1\text{ times}\}.\]

In particular, if $\omega\in \Omega$ then $\omega$ contains exactly two distinct symbols from $\{-1, 0, 1\}$ and one of those symbols occurs with multiplicity $1$. We, therefore, partition $\Omega$, by the symbol occurring exactly once, as $\Omega_{-1}\cup \Omega_{0}\cup \Omega_{1}$, where $\Omega_{\alpha}$ is defined for each $\alpha\in \{-1, 0, 1\}$ as \[\Omega_{\alpha}=\{\omega\in \Omega: \alpha \text{ occurs exactly once in }\omega\}\;.\]

For later use we make a convention to denote the elements of $\Omega_{\alpha}$. Let $i\in [n]$. Then, there are exactly two elements $\omega, \omega'\in \Omega_{\alpha}$ such that $\omega_{i}=\alpha=\omega_{i}'$. We arrange these elements in lexicographic order and call them $f_{\alpha, i}$ and $s_{\alpha, i}$ respectively. Note that $\Omega_{\alpha}=\cup_{i=1}^{n}\{f_{\alpha, i}, s_{\alpha, i}\}$. 

Let $\xi=(\xi_1, \ldots, \xi_n)$ be a random vector taking values in $\Omega$ and let $\mu$ denote the law of $\xi$. Notice that $\mu$ can be identified with a vector $\widetilde{\mu}\in \Delta_{6n}$ where $\Delta_{d}$ is the $d$-dimensional unit simplex. In other words, any collection of non-negative real numbers $p_{\alpha, i}\coloneqq \mu(f_{\alpha, i}), q_{\alpha, i}\coloneqq \mu(s_{\alpha, i})$ such that $\sum_{\alpha\in \{-1, 0, 1\}}\sum_{i=1}^{n}( p_{\alpha, i}+q_{\alpha, i})=1$ determines a unique random vector $\xi$ taking values in $\Omega$. 

We now record the following immediate fact.
\begin{fact}
    Let $\mu$ be the law of a random vector $\xi$ taking values in $\Omega$. Then, $\mu$ is exchangeable if and only if   $p_{\alpha, i}$ and $q_{\alpha, i}$ are independent of $i$
    for each $\alpha\in \{-1, 0, 1\}$ (that is $p_{\alpha, i}=p_{\alpha, j}$ for all $i, j\in [n]$ and each $\alpha$ and similarly for $q_{\alpha, \cdot}$).
\end{fact}

We now investigate the conditions for $\xi$ to be two-color exchangeable. To do this, let $f:\{-1, 0, 1\}\to \{0, 1\}$ be any non-constant function and let $\xi$ be a random vector taking value in $\Omega$. Then, the random vector $f(\xi)\coloneqq (f(\xi_1), \ldots, f(\xi_n))$ takes values in the set $\Sigma\subset \{0, 1\}^{n}$ where 
\[\Sigma = \left\{(x_1, \ldots, x_n): x_i\in \{0, 1\}, \sum_{i=1}^{n}x_i\in \{0, 1, n-1\}\right\}\;.\]

We further partition $\Sigma$ as $\Sigma_0\cup \Sigma_1\cup \Sigma_{n-1}$ where 
$\Sigma_T\coloneqq \{(x_1, \ldots, x_n)\in \Sigma: \sum_{i=1}^{n}x_i=T\}$ for $T\in \{0, 1, n-1\}$. Let $\eta$ be a random vector taking values in $\Sigma$ and let $\nu$ denote the law of $\eta$. Since each $\Sigma_T$ is an orbit under permutations of some binary string of length $n$, it follows that $\eta$ is exchangeable if and only if $\nu$ is constant on $\Sigma_T$ for each $T\in \{-1, 1\}$. Note that $|\Sigma_0|=1$ and therefore any $\nu$ is trivially constant on $\Sigma_0$. 

As argued in Section~\ref{sec:n=2}, $\xi$ is two-color exchangeable if and only if $f(\xi)$ to be exchangeable for each $f=\delta_{\alpha}$ with $\alpha=-1, 0, 1$. Let $\mu_{\alpha}$ denote the law of $f(\xi)$ where $f=\delta_{\alpha}$. Checking the condition that $\mu_{\alpha}$ is constant on $\Sigma_{T}$ for $T\in \{1, n-1\}$ yields the following six conditions:
\begin{enumerate}
    \item $\nu_{-1}$ is constant on $\Sigma_{1}$ if and only if $\mu(f_{0, i})+\mu(f_{1, i})$ is independent of $i$
    \item $\mu_{-1}$ is constant on $\Sigma_{n-1}$ if and only if $\mu(f_{0, i})+\mu(f_{1, i})$ is independent of $i$
     \item $\nu_{0}$ is constant on $\Sigma_{1}$ if and only if $\mu(f_{0, i})+\mu(s_{0, i})$ is independent of $i$
    \item $\mu_{0}$ is constant on $\Sigma_{n-1}$ if and only if $\mu(f_{-1, i})+\mu(s_{1, i})$ is independent of $i$
    \item $\nu_{1}$ is constant on $\Sigma_{1}$ if and only if $\mu(f_{1, i})+\mu(s_{1, i})$ is independent of $i$
    \item $\mu_{1}$ is constant on $\Sigma_{n-1}$ if and only if $\mu(s_{-1, i})+\mu(s_{0, i})$ is independent of $i$
\end{enumerate}

\subsection{A concrete example}
It is clear from the above discussion that there are two-color exchangeable sequences that are not exchangeable. For completeness, we specify a particular distribution $\mu$ on $\Omega$ that is two-color exchangeable but not exchangeable. 

To define $\mu$, it is enough to specify $\mu(f_{\alpha_i})$ and $\mu(s_{\alpha, i})$ for $\alpha\in \{-1, 0, 1\}$ and $i\in [n]$. To this end, we set 
\begin{align*}
    \mu(f_{-1, i}) = \mu(f_{1, i}) = \mu(s_{0, i}) &= \frac{n-i+1}{3n(n+1)}\\
    \mu(s_{-1, i}) = \mu(s_{1, i}) = \mu(f_{0, i}) &= \frac{i}{3n(n+1)}\;.
\end{align*}

The constant $3n(n+1)$ is normalizing constant to make $\mu$ a probability measure. Note that $\mu$ is not exchangeable because $\mu(f_{-1, i})$ depends on $i$, in particular, $\mu((-1, 0, \ldots, 0))\neq \mu(0, -1, \ldots, 0)$. However, notice that $\nu_{\alpha}(\omega)=\frac{1}{3n}$ for all $\omega\in \Sigma_1\cup \Sigma_{n-1}$ and for each $\alpha\in \{-1, 0, 1\}$. Thus, $\nu_{\alpha}$ is exchangeable for each $\alpha$. In particular, $\xi$ is two-color exchangeable. 

\section{Conclusion and further direction}
Constructing similar examples of infinite sequence that is two-color exchangeable but not exchangeable may be more delicate. We leave this investigation for future.

The notion of two-color exchangeability, even if it is not equivalent to exchangeablility, seems an interesting concept to study. Exchangeable sequences arise naturally in sampling without replacement. It would be interesting to find similar processes where two-color exchangeability arises naturally.

\bibliography{bibliography.bib}

\end{document}